\documentclass[11pt,reqno]{amsart}
\usepackage{amssymb,amsmath}\newcommand{\be}{\begin{eqnarray}}
\newcommand{\ee}{\end{eqnarray}}\newcommand{\e}{{\varepsilon}}\newcommand{\R}{{\mathbb R}}\newcommand{\Cant}{{\mathcal C}}\newcommand{\K}{{\mathcal K}}
\newcommand{\supp}{\operatorname{supp}}\newtheorem{theorem}{Theorem}\newtheorem{lemma}[theorem]{Lemma}\theoremstyle{definition}\theoremstyle{remark}\numberwithin{equation}{section}\input epsf.sty
\begin{document}\thispagestyle{empty}

%
%
%
%
%
%
%
%
%
%
%
%

\newcommand{\Anr}{A_{n,r}}
\newcommand{\Enr}{E_{n,r}}
\newcommand{\C}{\mathbb{C}}
\newcommand{\sigthet}{\sigma_{\theta}}
\newcommand{\sigthetf}{\sigthet ^f}
\newcommand{\sigthetg}{\sigthet ^g}
\newcommand{\sigtheth}{\sigthet ^h}
\newcommand{\Fasig}{Fav_\sigma}
\newcommand{\prosig}[1]{Proj_\theta(\sigma_\theta(#1))}
\newcommand{\proT}[1]{Proj_\theta(T_\theta(#1))}
\newcommand{\fnT}{f_{n,\theta,T}}
\newcommand{\rhoPT}{\rho_{P,T}}
\newcommand{\rhopsig}{\rho_{P,\sigma}}
\newcommand{\PP}{\mathcal{P}}
\newcommand{\Ent}{E_{n,\theta,T}}
\newcommand{\sigtil}{\tilde{\sigma} _\theta^\e}
\newcommand{\sigthetgtau}{\sigma_{\theta}^{g_{\tau}}}
\newcommand{\gtil}{\tilde{g}_{\tau _\theta}}
\newcommand{\rhopsigg}{\rho_{P,\sigma^g}}

\title[Intersections of large-radius circles with the four-corner Cantor set]{{Intersections of large-radius circles with the four-corner Cantor set: estimates from below of the Buffon noodle probability for undercooked noodles}}
\author{Matthew Bond}\address{Matthew Bond, Department of Mathematics, Michigan State University.
{\tt bondmatt@msu.edu}}
\author{Alexander Volberg}\address{Alexander Volberg, Department of  Mathematics, Michigan State University and the University of Edinburgh.
{\tt volberg@math.msu.edu}\,\,and\,\,{\tt a.volberg@ed.ac.uk}}

\begin{abstract}Let $\Cant_n$ be the $n$-th generation in  the construction of the middle-half Cantor set. The Cartesian square $\K_n$ of $\Cant_n$ consists of $4^n$ squares of side-length $4^{-n}$. The chance that a long needle thrown at random in the unit square will meet $\K_n$ is essentially the average length of the projections of $\K_n$, also known as the Favard length of $\K_n$. A result due to Bateman and Volberg \cite{BV} shows that a lower estimate for this Favard length is $c\,\frac{\log n}{n}$.

We may bend the needle at each stage, giving us what we will call a noodle, and ask whether the uniform lower estimate $c\,\frac{\log n}{n}$ still holds for these so-called Buffon noodle probabilities. If so, we call the sequence of noodles undercooked. We will define a few classes of noodles and prove that they are undercooked. In particular, we are interested in the case when the noodles are circular arcs of radius $r_n$. We will show that if $r_n \geq 4^{\frac{n}{5}}$, then the circular arcs are undercooked noodles.
\end{abstract}\maketitle

\section{{\bf Introduction}} \label{sec:intro}

Let $C_r(z):= \lbrace z + re^{i\theta}:\theta\in[0,2\pi]\rbrace$. We are interested in the Lebesgue plane measure of the set $\Anr := \lbrace z:C_r(z) \cap \K_n \not= \emptyset \rbrace $.
$$
|\Anr| =  \int_0^{2\pi}  \int_0^{\infty} \chi_{\Anr} (\rho e^{i \theta}) \rho\, d\rho d\theta \geq (r-2)  \int_0^{2\pi}  \int_{r-2}^{r+2} \chi_{\Anr} (\rho e^{i \theta}) d\rho d\theta.
$$

The last integrand above (excluding the $r-2$) can be thought of as a small translated distortion of the integrand of $Fav(\K_n)$. We will describe the distortion and show that for $r$ large enough, $\K_n$ is sufficiently coarse for the argument of Bateman-Volberg to yield the same lower bound.

To this end, let $T : \C\times S^1\rightarrow\C$, where we'll write for convenience $$T_\theta (z) := T(z,e^{i\theta})$$

Then define, for any $E\subset\C$, 
\begin{equation}\label{FaT def}Fav_T(E):=\frac{1}{2\pi}\int{\proT{E}}d\theta.\end{equation}
It remains to write $$\int \chi_{\Anr} (\rho e^{i\theta})d\rho = |\prosig{\K_n}|$$ for the appropriate choice of $\sigma$. Define \begin{equation}\label{f}f_r(y):= \begin{cases} r-\sqrt{r^2-y^2},\,|y| \leq 2\\r-\sqrt{r^2-4},\, \text{otherwise} \end{cases}\end{equation}

Then define $\sigma_0(x,y):=(x-f_r(y),y)$, and $\sigthet := R_{-\theta} \circ \sigma_0 \circ R_{\theta}$, where $R_\theta$ is clockwise rotation by the angle $\theta$. Then $(\rho + r) e^{i\theta}\in\Anr$ iff $\sigthet$ carries some point of $\K_n$ to the line perpendicular to $\theta$ at $\rho e^{i\theta}$, i.e., $\chi_{\Anr}((\rho + r) e^{i\theta}) = \chi_{\prosig{\K_n}}(\rho e^{i\theta})$. Thus \begin{equation}\label{Anr area} |\Anr| \geq 2\pi(r-2) \Fasig(\K_n) \end{equation}

Above, $f_r$ is an example of a $\textbf{noodle}$, which is a parameterized family of real functions. For any noodle $g$, we may define $\sigma_\theta^g$ from $g$ in the same manner that we defined $\sigthet$ from $f$. The symbol $\sigthet$ supresses $f$ and $r$ from the notation, but we will refer to them explicitly as needed.

\section{\bf Bateman-Volberg revisited} \label{Bateman-Volberg}

Let us review briefly the argument of Bateman and Volberg \cite{BV} which proves that $Fav(\K_n)\geq\frac{C\log n}{n}$. Below, we will fix an $n$, and none of the constants will depend on $n$. We rotate the axes, defining $\theta =0$ to be the direction $arctan(1/2)$, because $\K_n$ projects onto this direction nicely: the projected squares together fill out a single connected interval, and the projected squares intersect only on their endpoints. These almost-disjoint projected intervals induce a 4-adic structure on the interval.

For each square $Q$ of size $4^{-n}$ in $\K_n$, $\chi_{Q, \theta}(x) $ is the characteristic function of the projection onto the direction $\theta$. Put $f_{n, \theta}(x) =\sum_{Q, \ell(Q)=4^{-n}} \chi_{Q, \theta}(x)$. That is, $f_{n, \theta}(x)$ denotes the number of squares of length $4^{-n}$ whose
orthogonal projection on line $L_{\theta}$ contain a point $x$ of this line. Let us denote the support of $f_{n, \theta}(x)$  by $E_{n,
\theta}$, and let $|E_{n, \theta}|$ denote its length.

Let $J_j := (arctan(4^{-j}),arctan(4^{-j+1}))$. (The count starts from the special direction chosen above.) The central computation of Bateman-Volberg centers around a partitioning of an estimate of $Fav(\K_n)$ into conical neighborhoods $J_j\times \R$:

$$
\int_{J_j} |E_{n,\theta}|\, d\theta \ge\frac{(\int_{J_j}\int
f_{n,\theta}dx\,d\theta)^2}{\int_{J_j}\int
f^2_{n,\theta}dx\,d\theta}\
$$

Here we used the Cauchy inequality on $f_{n,\theta}$ and $\chi_{E_{n,\theta}}$.

Trivially, $\int_{J_j}\int f_{n,\theta}dx\,d\theta\leq C4^{-j}$. The interesting part of Bateman-Volberg amounts to showing that our partition has been chosen such that we may conclude that $\int_{J_j}\int f^2_{n,\theta}dx\,d\theta\leq Cn4^{-2j}$ for the approximately $\log n$ many values of $j$ ($3<j<\log n$), so that $\int_{J_j}{|E_{n,\theta}|d\theta}>C/n$, and summing over $j$ yields $Fav(\K_n)\geq\frac{C\log n}{n}$. Now

$$f_{n,\theta}^2=\sum_{Q,Q'}\chi_{Q,\theta}\chi_{Q',\theta}=\sum_{Q\neq Q'}\chi_{Q,\theta}\chi_{Q',\theta} + \sum_{Q}\chi_{Q,\theta}^2\,.$$

Integrating over $J_j\times\R$, the latter diagonal sum becomes $C4^{-j}\leq Cn4^{-2j}$ (the inequality uses $j<\log n$). When estimating the other integral, things become combinatorial - most of these terms are identically 0 in $J_j\times\R$. So define $A_{j,k}$ to be the set of pairs $P=(Q,Q')$ of Cantor squares such that in our special coordinate system, the centers $q$ and $q'$ of $Q$ and $Q'$ have vertical distance $4^{-k-1}\leq |y_q-y_{q'}| \leq 4^{-k}$ and satisfy the condition on horizontal spacing $4^{-j-1}\leq |\frac{x_p-x_{p'}}{y_p-y_{p'}}| \leq 4^{-j}$. We can think of $4^{-j}$ as being $tan(\theta)$ for $\theta$ such that the squares $Q,Q'$ overlap in the projection onto $\theta$. In Bateman-Volberg \cite{BV}, it was proved that

\begin{equation}
\label{Ajk}
|A_{j,k}|\leq C4^{2n-k-2j}
\end{equation}

For any $(j,k)$ pair, it is immediate that the integral $\rho_P:=\int_0^{2\pi}{\int_{\R}{\chi_{Q,\theta}\chi_{Q',\theta}d\theta dx}}$ satisfies $\rho_P\leq 4^{k-2n}$, and the integrand is supported only for angles belonging to $J_{j-1}$, $J_j$, and $J_{j+1}$. So we fix $j$ and sum over $k$ to get

$$\int_{J_j\times\R} {\sum_{Q\neq Q'}{\chi_{Q,\theta}\chi_{Q',\theta} d\theta dx}} \leq
$$
$$
 \sum_{k=1}^{n-j+1}{\max\lbrace\rho_P:P\in A_{j',k}\text{ for  } |j'-j|\leq 1\rbrace
(|A_{j-1,k}| +|A_{j,k}|+|A_{j+1,k}|)}\leq Cn4^{-2j}.$$

(Note above that $j+k\leq n$, since $4^{-j-k}$ bounds the horizontal distance between centers of squares from above.)

This completes the proof of the result of Bateman and Volberg. We will need to remember some of the notations for later, and the estimate \eqref{Ajk}.

\section{\bf A simple lemma} \label{angle lemma}

Now we show that the $\sigthet$ in the integrand of $\Fasig(\K_n)$ hardly disturbs the angular sorting argument of Bateman-Volberg. We will need the following estimate on $|f_r'(y)|$:
\begin{equation}\label{f'} |f_r'(y)|\leq\frac{4}{r}\end{equation}
because it gives us
\begin{equation}\label{lipf} Lip(\sigthet-Id)\leq\frac{4}{r}\end{equation}
when we conjugate $\sigma_0$ by the isometry $R_\theta$.

\begin{lemma}
\label{angles} Let $\e >0$ be small enough. Let $T:\C\to\C$ be such that $Lip(T-Id)<\e$. Then $\forall z,w\in\C$, $$|arg(z-w)-arg(T(z)-T(w))|<2\e $$ (for appropiate choices of arg).
\begin{proof}
Write $z-w=\rho e^{i\theta}$, and let $\alpha:=arg(z-w)-arg(T(z)-T(w))$. $$arg(T(z)-T(w))=arg((T-Id)(z)-(T-Id)(w)+(z-w))=arg(\lambda\rho e^{i\beta}+\rho e^{i\theta})$$ for some $\lambda<\e,\beta\in[0,2\pi]$. So $arg(T(z)-T(w))=arg(\lambda e^{i\beta}+ e^{i\theta})$

Then $|\alpha | \leq \hat{\alpha}$, where $tan(\hat{\alpha})=\frac{\e}{1-\e}\Rightarrow |\alpha|<2\e$.
\end{proof}
\end{lemma}

\section{A few classes of undercooked noodles}\label{theorems}

We say that $T_n:\C\times S^1\to\C$ is an $\textbf{undercooking}$ of the plane if $Fav_{T_n}(\K_n)\geq C\frac{\log n}{n}$. Likewise, we say that $\lbrace r_n \rbrace$ is $\textbf{undercooked}$ if $\sigma^{f_{r_n}}$ is an undercooking of the plane. In fact, this is the same as saying that $f_{r_n}$ is an undercooked noodle.

\begin{theorem}
\label{4^n/5}
If $r_n \geq 4^{n/5}$, then $r_n$ is an undercooked sequence.
\end{theorem}

First we will prove a more general result which is weaker in the sense that it does not give us the above theorem unless we strengthen the $4^{n/5}$ in the hypothesis of the above theorem to $4^n$.

\begin{theorem}
\label{4^n}
If $T_n:\C\times S ^1\to\C$ satisfies $Lip(T_{n,\theta}-Id)<4^{-n} \,\,\forall n,\theta$, then $T_n$ is an undercooking of the plane.
\end{theorem}

Note that $T_n$ need not be induced by a noodle.

\section{A sorting lemma and the weak $\rho_P$ estimate}\label{FaT lemma}

For any $T:\C\times S^1\to\C$, define $A_{j,k,T}$ by $P=(Q,Q')\in A_{j,k,T}$ if and only if $\exists\theta : (T_\theta (Q), T_\theta (Q'))\in A_{j,k}$

\begin{lemma}
\label{AjkT}
$\textbf{Sorting Lemma}$ \\
Let T satisfy $Lip(T_\theta-Id)<\frac{1}{8n}$. Then $\forall j<\log n$, $|A_{j,k,T}|\leq C4^{2n-k-2j}$.
\begin{proof}
Distances are preserved up to a multiple of $1\pm \frac{1}{n}$ under $T$, so for a $j,k$ pair, $k$ can change by at most one under $T_\theta$. Lemma \ref{angles} implies that angles are changed additively by at most $\frac{1}{4n}$ under $\sigma$, so $j$ can change by at most one if $j\leq \log n$. Thus Bateman-Volberg \eqref{Ajk} gives us 
$$|A_{j,k,T}|\leq \sum_{-1\leq l,m \leq 1}|A_{j+l,k+m}|\leq C4^{2n-k-2j}.$$
\end{proof}
\end{lemma}

Note that $T=\sigma^f$ satisfies Lemma \eqref{AjkT} for $r>32n$, but this will NOT be sufficient for the  $\rho_P$ estimate.

Instead of $f_{n,\theta}$ and $\rho_P$, consider

$\fnT := \sum_Q{\chi_{T_\theta (Q),\theta}}$ and $\rhoPT := \int|\proT{Q} \cap \proT{Q'} |d\theta$

\begin{lemma}
\label{weak rhoP}
$\textbf{Weak $\rho_P$ Lemma}$ \\
Let T be as in Theorem \ref{4^n}. Then $\rhoPT\leq 4^{k-2n}$
\begin{proof}
$T$ at most stretches by $1+4^{-n}$. We write $4^{-n}=\frac{4}{r}$ both as an abstraction and to anticipate Theorem \ref{4^n/5}.

It is immediate that for two squares of size $4^{-n}$ at distance $\asymp 4^{-k}$ one has $|\lbrace\theta :Proj_\theta (Q)\cap Proj_\theta (Q')\rbrace |\leq C4^{k-n}$, so Lemma \ref{angles} implies 
\begin{equation}\label{*}|\lbrace\theta :\proT {Q}\cap\proT{Q'}\rbrace |\leq C(4^{k-n} +1/r)\end{equation}
and here we use $r\geq 4^n$ to conclude $|\lbrace\theta :\proT {Q}\cap\proT{Q'}\rbrace |\leq C4^{k-n}$, and thus the lemma as the length of projections is obviously bounded by $C\,4^{-n}$ .
\end{proof}
\end{lemma}

\section{Proof Theorems \ref{4^n/5} and \ref{4^n}}
\label{ProofofTh}

Theorems \ref{4^n/5} and \ref{4^n} can now be proved in the spirit of Bateman-Volberg. However, Week $\rho_p$ lemma is much too weak for Theorem \ref{4^n/5}, an analogous strong $\rho_P$ lemma will be needed in the case of Theorem \ref{4^n/5}. We will state that lemma now and prove it later.

\begin{lemma}
\label{strong rhoP}
$\textbf{Strong $\rho_P$ Lemma}$ \\
Let $r_n\geq 4^{n/5}$ (as in Theorem \ref{4^n/5}0. Then $\rhopsig\leq 4^{k-2n}$.
\end{lemma}

\medskip

Take this lemma for granted to finish the proof of Theorem \ref{4^n/5}.

\medskip

Let $\PP _{j,T}:=\bigcup^{n-j}_{k=0}A_{j,k,T}$. Then $\sum_{P\in{\PP}_{j,T}}{\rhoPT}\leq Cn4^{-2j}$ (Sorting and $\rho_P$ Lemmas).
Also, let $\Ent := \supp \fnT$. A couple applications of the Cauchy inequality to $\fnT$ and $\chi_{\Ent}$ give us 
\begin{equation}
\label{Jj int}
\int_{J_j}|\Ent |\geq \frac{(\int_{Jj}\int_\R\fnT dxd\theta)^2}{(\int_{Jj}\int_\R\fnT ^2dxd\theta)}
\end{equation}

We have \begin{equation}
\label{fnj int} \int_{Jj}\int_\R\fnT dxd\theta \approx 4^{-j}
\end{equation}
$$\int_{Jj}\int_\R\fnT ^2 dxd\theta = \int_{Jj}\int_\R\fnT dxd\theta + \sum_{Q\neq Q'}{\int_{Jj}\int_\R\chi _{T_\theta (Q)} \chi _{T_\theta (Q')} dxd\theta}$$
$$\leq C4^{-j} + \sum_{P\in{\PP}_{j-1}\cup {\PP}_{j}\cup {\PP}_{j+1}}{\rho_{P,T}}\leq C(4^{-j} + n4^{-2j}) \leq C\,n\,4^{-2j}\,,$$
where the last inequality relies on $j<\log n$. So $(\ref{Jj int})$, $(\ref{fnj int})$ give us, together with the above, $\int_{Jj}{|\Ent |d\theta} \geq C/n.$
Summing over $3<j<\log n$, we get the result.
$\square$

\section{Some useful facts about shear group}
\label{sheargroup}

We need to prove the Strong $\rho_P$ Lemma. Before we proceed, a few facts about shear groups need to be stated. Below, $g$ and $h$ will be arbitrary noodles. Recall that $\sigma_0^g(x,y):=(x-g(y),y)$, and $\sigthetg := R_{-\theta} \circ \sigma _0^g\circ R_\theta$.
First, there is this simple fact for arbitrary functions $g$ and $h$:

$$\sigthetg\circ\sigtheth = \sigthet ^{g+h}$$
Next, we show how shears by linear noodles behave.
For $g(y)=b$, we get

\begin{equation}
\supp(Proj_\theta (\sigthetg (E_\theta)))=\supp(Proj_\theta (E_\theta))-b
\end{equation}

For $g(y)=my$, $\alpha:= \arctan m$, we get 
\begin{equation}
\supp(Proj_\theta (\sigthetg (E_\theta)))=\bigr(R_{\alpha}\frac{\supp(Proj_{\theta -\alpha}(E_\theta))}{cos(\alpha)}\bigr)=(\sqrt{1+m^2})R_{\alpha} \supp(Proj_{\theta - \alpha}(E_\theta))
\end{equation}

For for $g(y)=my +b$, then, given a set $A$ on the real line,
\begin{equation}
\int_0^{2\pi}\int_A\chi_{Proj_\theta (\sigthetg (E_\theta))}(x)dxd\theta=\sqrt{1+m^2}\int_0^{2\pi}\int_{\frac{1}{\sqrt{1+m^2}}(A+b)}\chi_{Proj_\theta (E_{\theta-\alpha})}dxd\theta 
\end{equation}

\section{Proof of the Strong $\rho _P$ Lemma}
\label{strong}

Recall: $f(y)=r-\sqrt{r^2-y^2}$, and $|f'(y)|<C/r$.
We also have $f''(y)=\frac{r^2}{(r^2-y^2)^{3/2}}$, so $|f''(y)|<C/r$. Remember that $\sigthet$ refers to $\sigthetf$ if no noodle is specified. \\
Remember that we still have \eqref{*}: $$|\lbrace\theta :\sigthetf {Q}\cap\sigthetf{Q'}\rbrace |\leq C(4^{k-n} +1/r).$$
$r\geq 4^{n/5}$, so we are done proving that this measure is bounded by $C\,4^{k-n}$ for all $k\geq 4n/5$. So let $k<4n/5$.\\
WLOG, the centers of $Q$ and $Q'$ are $(0,0)$ and $(0,-L)$. To see this, note that $$\rhopsig\approx \frac{1}{r}\int_ 0^{2\pi}\int \chi _{\lbrace Q\cap C_r (\rho e^{i\theta})\neq\emptyset\rbrace} \chi _{\lbrace Q'\cap C_r(\rho e^{i\theta})\neq\emptyset\rbrace} \rho d\rho d\theta$$
$$=\frac{1}{r}\int_ 0^{2\pi}\int \chi _{\lbrace Q^*\cap C_r (\rho e^{i\theta})\neq\emptyset\rbrace} \chi _{\lbrace Q'^*\cap C_r(\rho e^{i\theta})\neq\emptyset\rbrace} \rho d\rho d\theta \approx \rho _{P^*,\sigma},$$
where $P^*$ is the pair $(Q,Q')$ translated and rotated to $(Q^*,Q'^*)$ as in the WLOG condition.
(The area of the set of centers of circles for which the indicated intersections occur is obviously invariant under translations and rotations of the plane, and the possible $\rho$-values for which the intersection occurs are restricted to an annulus of inner and outer radius $\approx r$. Thus, the $\rho$ in $\rho d\rho d\theta$ is $\approx r$, both before and after the translation and rotation described.) \\

As $\theta$ ranges over all angles such that $Q,Q'$ have intersecting $\sigthet$-projections, 
the angle distortion of Lemma \ref{angles} says that such angles $\theta$ satisfy $|\theta|< \frac{C}{r} + C, 4^{k-n}<\frac{C}{r}$ (see \ref{*}) as $k<4n/5$.

For these $\theta$, rotation $R_{\theta}(Q)$ is in the band $\delta \leq y \leq L + \delta $, for $\delta = 4^{-n} + L(1-cos(C/r))$, giving $\delta \leq C\max\lbrace 4^{-n}, L/{r^2}\rbrace \leq C4^{-2/{5n}}$. Transform the integral using the shear group.
Let $l(y)$ linearly approximate $f(y)$ at $y=L-\delta$, with $l(y)=my+b$. Note that $|b|\leq CL/r$. Let $\e (y):=f(y)-l(y)$ on $[L-\delta , L+\delta]$ and extend $\e$ continuously to be constant elsewhere. Then, with $b':=b/\sqrt{1+m^2}$:

$$\rhopsig =\int{|Proj_\theta (\sigthetf (Q')) Proj_\theta(\sigthetf (Q))|d\theta}\leq \int_0^{2\pi}\!\!\int_{-4^{-n}}^{4^{-n}}\chi_{Proj_\theta (\sigthetf (Q'))}(x)\,dxd\theta$$
$$=\int_0^{2\pi}\!\!\int_{l_{\theta}\cap [-4^{-n}}^{4^{-n}]}\chi_{Proj_\theta (\sigthet ^l (\sigthet ^\e (Q')))}\,dxd\theta\leq C\,\int_0^{2\pi}\!\!\int_{l_{\theta}\cap [b'-2\cdot 4^{-n}, b'+2\cdot 4^{-n}]}\chi_{R_{\alpha}Proj_{\theta-\alpha} (\sigma_{\theta} ^\e (Q'))}\,dxd\theta\,.$$
Changing variable, we see that this is at most
$$
C\,\int_0^{2\pi}\!\!\int_{l_{\theta+\alpha}\cap [b'-2\cdot 4^{-n}, b'+2\cdot 4^{-n}]}\chi_{R_{\alpha}Proj_{\theta} (\sigma_{\theta+\alpha} ^\e (Q'))}\,dxd\theta\,.
$$

Let $\Gamma := \lbrace\theta : R_{\alpha}Proj_{\theta} (\sigma_{\theta+\alpha} ^\e (Q'))\cap l_{\theta+\alpha}\cap [b'-2\cdot 4^{-n},b'+2\cdot 4{-n}]\neq\emptyset\rbrace$, and let $z:=(0,-L)$.
If $\theta\in\Gamma$, then $R_{\alpha}Proj_{\theta} (\sigma_{\theta+\alpha} ^\e (z))\in  l_{\theta+\alpha}\cap [b'-3\cdot 4^{-n}, b'+3\cdot 4^{-n}]$.

Using $|f''(y)|<C/r$, we get $|\e '(y)|<C\delta /r < C\,L/r^3<C4^{-3/{5n}}$. Then it follows that $|\e (y)|<C\,\delta^2/r < C\, L^2/r^5 < C4^{-n}$. So $|\sigma_{\theta'+\alpha}^\e (z) - z|< c\, 4^{-n}$, and hence $|R_{\alpha}Proj_{\theta} (\sigma_{\theta+\alpha}^\e(z)) - R_{\alpha}Proj_{\theta} (z)| \leq C\,4^{-n} \,\,\forall \theta\in\Gamma$. So 

$$\Gamma\subseteq\lbrace\theta : Proj_{\theta} (z)\in  R_{-\alpha}(l_{\theta+\alpha}\cap [b'-C4^{-n}, b'+C4^{-n}])\rbrace =
$$
$$
\lbrace\theta' : L\sin\theta\in[b'-C4^{-n},b'+C4^{-n}]\rbrace\,,
$$
which implies:
\begin{equation}\label{Gamma}|\Gamma |\leq C|\lbrace\theta : \sin\theta\in [b/L -C4^{k-n}, b/L +C4^{k-n}]\rbrace|.\end{equation}
Since $b << L$ and $k<4n/5$, $\sin\theta\approx\theta$, and we get $|\Gamma |\leq C4^{k-n}$, completing the proof of the Strong $\rho_P$ Lemma.
$\square$

\section{General Buffon noodle probabilities and the $\rho _P$ lemmas for arbitrary noodles}
\label{general noodles}

Let us define general noodle probabilities now. Let $g_{\tau} (y) := g(y-\tau)$. For a probability distribution $P$ on $\R ^2\times S^1$, a set $E\subset\C$, and noodle $g$, we can define $$Bu^g(E)=\int{proj_{\theta} (\sigthetgtau (E))(x)dP(x,\tau ,\theta)}.$$ We can choose an $L>10$, say, and let $P$ be normalized Lebesgue measure on $(-2,2)\times (-L,L) \times (0,2\pi )$, under which $$Bu^g(E) = \frac {1}{16\pi L}\int_0^{2\pi} {\int_{-L}^L {|Proj_\theta (\sigthetgtau (E))|d\tau d\theta}}= \frac{1}{16\pi L}\int_{-L}^L{Fav_{\sigthetgtau} (E) d\tau}.$$

Having done this, we will say that a noodle $g_n$ is $\textbf{undercooked}$ if $Bu^{g_n}(\K_n)\geq C \frac{\log n}{n}$.

Next, we describe the portion of the domain of integration in which the noodle hits the center of a square $Q$ at the same point $-\tau _0$ of the noodle. That is, if $Q$ has center $z=\rho e^{i\theta_0}$, consider $\tilde{g} := g-g(-\tau _0)$ and $\sigma _\theta ^{\tilde{g}_{\tau _0}}$. For each $\theta$, we need to find the unique $x _\theta$ and $\tau _\theta$ such that the line centered at $x_\theta e^{i\theta}$ and with positive axis in the $\theta + \pi /2$ direction intersects $z$ at $y=\tau _\theta - \tau _0$. In fact, $x_\theta =|z|cos(\theta -\theta_0)$ and $\tau_\theta= \tau_0 - |z|\sin(\theta -\theta_0)$.
(Diagram) \\ Then when computing

$$\int_0^{2\pi} {\int_{x_\theta -a}^{x_\theta + a} {Proj_\theta ( \sigthet ^{\gtil} (E))(x) dxd\theta}},$$

WLOG $z=0$. That is, 

\begin{equation}
\label{noodle WLOG}
\int_0^{2\pi} {\int_{x_\theta -a}^{x_\theta + a} {Proj_\theta ( \sigthet ^{\gtil} (E))(x) dxd\theta}} = \int_0^{2\pi} {\int_{-a}^{a} {Proj_\theta ( \sigthet ^{\tilde{g}} (E-z))(x) dxd\theta}}.
\end{equation}

So define $\rhopsigg =\int_{-L}^L{\int_0^{2\pi}{|Proj_\theta (\sigthetgtau (Q)) Proj_\theta (\sigthetgtau (Q'))|d\theta d\tau}}$. We want $\rhopsigg < C4^{k-2n}$.

For $z =$ center of $Q$, and for fixed $\tau_0$, define $D = \lbrace \tau = \tau_0 -|z|\sin(\theta - \theta_0), |x-|z|cos(\theta - \theta_0)|\leq C4^{-n}, \theta\in (0,2\pi) \rbrace$. Then if $I_D(\tau_0) := \int_D{Proj_\theta (\sigthet ^{\gtil}(Q'))(x)dxd\theta}$, then $\rhopsigg \leq \int_{-L}^{L}{I_D(\tau_0) d\tau_0}$. So because of \eqref{noodle WLOG}, we are in the same case as the Strong $\rho_P$ Lemma for circles, so long as the estimates $|g_n(y)|<1$, $|g_n'(y)|<4^{-n/5}$, and $|g_n''(y)|<4^{-n/5}$ hold. Likewise, the Weak $\rho_P$ Lemma generalizes here so long as $|g_n(y)|<1$ and $|g_n'(y)|<4^{-n}$. In particular, in either case, such $g_n$ are undercooked.

A careful examination of the Strong $\rho _P$ Lemma shows that we can be slightly more flexible, requiring that the quantity $||g_n'(y)||_\infty ^4 ||g_n''(y)||_\infty <4^{-n}$ and $|g_n'(y)|<1/100$, instead. Using this, we get, for example, the undercooked noodle $g_n(y)=4^{-n/2}sin(4^{n/4}y)$.

%
%
%
%

\section{Closing remarks}
\label{remarks}

The above arguments are very local in nature, and fail to allow any large-scale bending. It is currently unclear what other sequences may or may not be undercooked - even for constant sequences, it is unclear. Since the random Cantor sets of \cite{PS} decay in Favard length like $C/n$ almost surely, perhaps if $r$ is small compared to $n$, this "randomizes" the Cantor set to the point of making $C/n$ an upper bound.

  \bibliographystyle{amsplain}

\end{document}